\documentclass[12pt]{amsart}

\def\call#1{\ensuremath{{\mathcal #1}}}
\def\M#1{\ensuremath{\mathbb #1}}
\def\f{\ensuremath{\varphi}}
\def\wt#1{\ensuremath{\widetilde{#1}}}

\def\cyl#1{\ensuremath{\textrm{Cyl}(#1)}}

\newtheorem{teo}{Theorem}[section]
\newtheorem{prop}[teo]{Proposition}
\newtheorem{cor}[teo]{Corollary}
\newtheorem{lemma}[teo]{Lemma}
\newtheorem{defi}[teo]{Definition}
\newtheorem{remark}[teo]{Remark}

\newcounter{step}
\def\thestep{\arabic{step}}
\newcommand{\step}{\refstepcounter{step}\noindent{\em Step \thestep.~}}
0

\begin{document}
\title[Length compactness for automorphisms of free groups]
{Geodesic currents and length compactness for automorphisms of
  free groups} 

\author{Stefano Francaviglia}
\thanks{The author was supported by a Marie Curie Intra European Fellowship.}
\thanks{
  currents}
\address{Departament de Matem\`atiques\\ Edifici C\\
Universitat Aut\`onoma\\
08193 Bellaterra
(Barcelona) Spain.}
\email{s.francaviglia@sns.it}

\begin{abstract}
We prove a compactness theorem for automorphisms of free
groups. Namely, we show that the set of automorphisms keeping
the length of the uniform current bounded is compact (up to conjugations.) This
implies that the spectrum of the length of the images of the uniform
current is discrete, proving a conjecture of I. Kapovich. 
\end{abstract}

\maketitle

\tableofcontents
\section{Introduction}

The aim of this paper is to investigate the action of the
automorphisms of a free group $F$ on the Cayley graph of $F$. In
particular, we are interested in understanding how automorphisms can
stretch geodesics. One can define the {\em length} of an automorphism as the
generic stretching factor (see~\cite{kks06}) which is, roughly
speaking, the average of the stretching ratios, taken over all
geodesics (Definition~\ref{defL}.) Intuitively, 
for $\Phi\in{\rm Aut}(F)$, one can think of its length as the
limit, as $n$ tends to infinity, of $||\Phi(w_n)||/n$ where $w_n$ is a
``random'' cyclically reduced word in $F$ of length $n$, and
$||\cdot||$ denotes the cyclically reduced length.

The length function on Aut$(F)$ is invariant under conjugation, so it
descends to a length function on Out$(F)$. Our main
result is the following:

\begin{teo}[Length compactness theorem]\label{t2}
  Let $\Phi_n$ be a sequence of automorphisms of $F$. Then, up to passing
  to subsequences, there exists a sequence $v_n\in F$ such that the
  automorphisms $\Psi_n$ defined by $x\mapsto v_n\Phi(x)v_n^{-1}$
  satisfy one (and only one) of the following:
  \begin{itemize}
  \item $\Psi_n$ converges to an automorphism $\Phi$ (that is, it has a
  constant subsequence $\Psi_{n_i}=\Phi$.)
  \item $L(\Phi_n)$ goes to $\infty$.
  \end{itemize}
\end{teo}

A reformulation of this theorem is the following:

\vskip\baselineskip

\noindent{\em The set of automorphisms of bounded length is compact
  up to conjugations. That is, for any $M\in\M R$, the set $\{[\Phi]\in
  $Out$(F):L(\Phi)<M\}$ is finite.} 

\vskip\baselineskip
Equivalently:
{\em For a sequence of automorphisms $\Phi_n$, if there is a
  word $w$ such that the cyclically reduced length of $\Phi_n(w)$ goes
  to $\infty$, then $L(\Phi_n)\to\infty$.}

\vskip\baselineskip

All the work pivots on the fact that the Cayley graph of $F$ is a
hyperbolic object. Therefore, the boundary at infinity $\partial F$ of
$F$ is well defined, and encodes enough information about the
dynamic of the action of Aut$(F)$. The main idea is that
length controls attractors: if $\Phi_n$ is a sequence of automorphisms
of bounded length then, up to conjugation, there are no attractors for
the action of $\Phi_n$ on $\partial F$. Using this fact we prove
that the sequence $\Phi_n$ keeps the cyclically reduced 
length of any element of $F$ bounded, and this implies that $\Phi_n$ has a
constant subsequence.    

As a consequence of Theorem~\ref{t2}, we get the following result. 
\begin{cor}\label{c2.1_2.2}
The spectrum of the length function is discrete. That is, the set 
$$\{L(\Phi):\Phi\in\textrm{Aut}(F)\}$$ 
is a discrete subset of
$\M R$.
\end{cor}

We remark that, while the length of an automorphism depends on the
free basis chosen for the Cayley graph, the spectrum of the length
function is intrinsic (it depends in fact on the current we use to
define it, which in our case is the uniform current.)

Corollary~\ref{c2.1_2.2} was conjectured to be true by I. Kapovich,
 inspired by computational evidence and partial results. For example,
 in~\cite{kks06} V. Kaimanovich, I. Kapovich and P. Shupp proved 
 (among other results) that an automorphism of length one must be simple
 (see below) and estimated  the ``first gap'' of the length function.

A consequence of Corollary~\ref{c2.1_2.2} is that one can use
inductive arguments on the length. For example, we get 
the following result, which
can be viewed as an {\em Ideal Whitehead Algorithm}
(see~\cite[Conjecture 5.3]{kapre05}.)
Recall that, given a free basis $\Sigma$ of $F$, an automorphism
$\tau$ is simple
if it is either a permutation of $\Sigma$ or an inner automorphism,
while it is a Whitehead automorphism of second kind if there is
$a\in\Sigma$ such that 
$\tau(x)\in \{x,xa,a^{-1}x,a^{-1}xa\}$ for all $x\in\Sigma$ (see for
example~\cite{lysc:libro} for more details.)

\begin{teo}\label{t_2.4_1.3}
  Let $\Phi\in\textrm{Aut}(F)$ be a non-simple automorphism. Then
  there exists a factorisation
$$\Phi=\tau_n\tau_{n-1}\cdots\tau_1\sigma$$
where $n\geq1$, the automorphism $\sigma$ is simple, each $\tau_i$ is
a Whitehead automorphism of the second kind, and
$$L(\tau_{i-1}\cdots\sigma)<L(\tau_i\tau_{i-1}\cdots\sigma)
\qquad i=1,\cdots,n-1.$$
\end{teo}

Let us say a few words about Theorem~\ref{t_2.4_1.3}. The {\em automorphism
  problem} for a free group $F$ asks, given two arbitrary elements
  $u,v\in F$, whether there exists $\Phi\in$ Aut$(F)$ such that
  $\Phi(u)=v$. In~\cite{whi36} Whitehead gave an algorithm solving
  that problem. The first part of the algorithm is to maximally 
  reduce the lengths of
  $u$ and $v$ via Whitehead automorphisms. Then, given two
  {\em minimal} elements one proves that they are in the same
  Aut$(F)$-orbit if and only if they are related via a sequence of minimal
  elements, each one obtained from the preceding via a Whitehead
  automorphism. Roughly speaking, Theorem~\ref{t_2.4_1.3} is an {\em
  averaged} version of the first part of the Whitehead algorithm.  We
  refer the reader to~\cite{kapre05} for a more detailed discussion on
  the matter. We only notice that, as our proof is ``typically
  hyperbolic,'' one may expect that it could be adapted to a more general
  setting, like that of hyperbolic
  groups, for which the automorphism problem is still not completely solved.

\vskip\baselineskip

The main tool we use is the theory of {\em geodesic currents}. These
are locally finite $F$-invariant Borel measures on the space of
geodesics lines in the Cayley graph of $F$.
Geodesic currents where introduced by F. Bonahon~\cite{Bon86}
in the setting of hyperbolic manifolds, and turned out to be very useful
in group theory (see for
example~\cite{Mar:tesi,kapre04,kapre05}). 

We also consider measures on the space of geodesic rays (that is,
half geodesics) in the Cayley graph of $F$. 
This is the space of {\em frequency measures} (see also~\cite{kap05}.)

Such spaces are in fact homeomorphic, but each one
has peculiar characteristics that are well-adapted to different situations, 
and we  
will jump from one to the other depending on the calculations we will
need to do. Roughly speaking, the action of Aut$(F)$ on currents is
``more natural'', while frequency measures are ``more compact''.

\vskip\baselineskip
\noindent As mentioned earlier, the length function on Aut$(F)$ depends on the
choice of a current. We use the uniform current, which is the
analogous of the Liouville current for the geodesic flow of a
surface but, {\em \c ca va sans dire}, length compactness should hold for many
other currents. In the last section we discuss some generalisations of
Theorem~\ref{t2} in such a direction.

\begin{remark}
The
space of automorphisms of a free group is discrete. Thus, compactness
is equivalent to finiteness, and to say that a sequence converges is
equivalent to 
say that it is finite (and hence eventually constant.) Nevertheless, we
prefer to speak about compactness and convergence because we think
that this is closer to the spirit of the paper, in which
we used ``more hyperbolicity than  discreteness'' (Even if
discreteness is necessary, as, for example, in Corollary~\ref{c2.1_2.2}.)
\end{remark}

\vskip\baselineskip
\noindent\textsc{Acknowledgement.}
I warmly thank Pepe Burillo, Bertrand Deroin, Warren Dicks, 
Jean-Fran\c cois Lafont, Joan
Porti, Enric
Ventura, Asli Yaman and all the guys of the {\em geometric group
  theory research group} of the UAB. I am in debt with
Armando Martino for his many comments and
suggestions. 

It is a pleasure to thank Ilya  Kapovich, for his useful observations
and for having pointed out a gap in a previous version of the
present work.

\section{Notation}\label{s_n}

For the remainder of the paper, we fix the following notation:
\begin{itemize}
\item $F$ is a free group of rank $k$, with a fixed free basis $\Sigma$. 
We set $A=\Sigma\cup \Sigma^{-1}$. Any element of $F$ corresponds to a unique
{\em freely reduced} word in the alphabet $A$, that is, a word not
containing sub-words of the form $aa^{-1}$ with $a\in A$. 
We identify $F$ with the
set of freely reduced words. We denote by $1$ the neutral element of $F$
(the empty word.)  A word $w$ is {\em cyclically reduced} if all the cyclic
permutations of $w$ are freely reduced.
 For $v$ a freely reduced word, $|v|$ denotes its length, and $||v||$
 denotes its {\em cyclically reduced length} that is, the length of
 the cyclically reduced word obtained by 
 cyclically reducing $v$. 
\item The Cayley graph of $F$ corresponding to $A$ will be shortly
  denoted simply by {\em Cayley graph}. We denote 
  by $1$ the base point of the Cayley graph corresponding to $1\in F$.
\item $\partial F$ is the {\em boundary at infinity} of $F$, identified with
 the set of geodesic rays in the Cayley graph, that is, freely
 reduced, right-infinite 
 words in the alphabet $A$. The boundary $\partial F$ is endowed with
 the Cantor-set topology. Namely, for each $v\in F$ we denote by $\cyl
 v$ the set of rays having $v$ as initial segment. We set
 $\cyl1=\partial F$. Then, a basis for the topology of $\partial F$ is
 given by the sets $\{\cyl v:v\in F\}$.
\item $\partial^2F$ is the set $\{(x,y)\in(\partial F)^2:x\neq
  y\}$. We identify $\partial^2F$ with the set of oriented bi-infinite
  geodesics in the 
  Cayley graph. We define the {\em base-ball} $B$ of $\partial^2F$ 
  as the set of geodesics passing through $1$.
\item For any $x\neq y\in F$, the cylinder $\cyl{[x,y]}$ is defined as
  the subset of $\partial^2 F$ of geodesics passing through the
  oriented segment joining $x$ and $y$ in the
  Cayley graph (with the correct orientation.) We set $\cyl{[1,1]}=B$.

\item We denote by $T:\partial F\to \partial F$ the shift operator deleting the
  first letter of a ray. It turns out that $T$ is a continuous map.
\item Given a topological space $M$, we identify the space of Borel
  measures on $M$ with the dual of $\call{C}_0(M)$ (the space of continuous
  functions on $M$ with compact support) endowed with the weak-*
  topology. Namely, measures $\mu_i$ converge to $\mu$ if and only if
  $\int \f\, d\mu_i\to\int\f\, d\mu$ for all $\f\in\call{C}_0(M)$. 
  If $\mu$ is a Borel measure on $M$, $N$ is a topological
  space, and $f:M\to N$ is a measurable, proper map, we denote by
  $f_*\mu$ the push-forward of $\mu$, that is the measure on $N$ such
  that $\int_N\f\,d(f_*\mu)=\int_M\f\circ f\,d\mu$ for all $\f\in\call
  C_0(N)$ (see for example the first chapters of~\cite{fed:libro}
  or~\cite{afp:libro} for good introductions to geometric 
  measure theory.)   

\end{itemize}

\section{Definitions and preliminary facts}\label{s3}

In this section, we define the space of geodesic currents and
of frequency measures, and we show that such spaces are homeomorphic.
We introduce the {\em length} of a current, which is the
analog of the length of a cyclically reduced word. 
We define the {\em uniform frequency measure} and the {\em
  uniform current}, which we use to define the  {\em length} of
  automorphisms.

First of all,  in order to describe the action of Aut$(F)$ on
currents, we need the 
following classical result, whose proof can be found in~\cite{coo87}. 
\begin{teo}
Let $\Phi$ be an automorphism of $F$. Then it  extends
 to a homeomorphism of $\partial F$ (which we still denote by $\Phi$.)
\end{teo}
Since $\Phi$ is a homeomorphism of $\partial F$, the map $\Phi\times\Phi$
is continuous and proper on $\partial^2F$.
It follows that any automorphism $\Phi$ acts on the space of Borel
measures on $\partial^2F$ via
$(\Phi\times\Phi)_*$.
The inclusion of $F$ in Aut$(F)$ given by inner automorphisms induces
an action of $F$ on the space of Borel measures on $\partial^2F$.
By abuse of notation, if $\eta$ is a Borel measure,
 we will denote by $\Phi\eta$ the measure
$(\Phi\times\Phi)_*\eta$. We can now give the definition of currents
and frequency measures. Our definitions are a little different from
those introduced in~\cite{kapre04,kap05}, as we do not require
measures to be normalised to probability measures. This is because the
quantities we are interested in (lengths of automorphisms) depend
on the total mass of the measures we work with.

\begin{defi}[Geodesic currents and their lengths]
  The space of {\em geodesic currents} is the space of
  locally finite ({\em i.e.} finite on compact sets) $F$-invariant
  non-negative Borel measures on $\partial^2 F$. 
 The {\em length} $L(\eta)$ of a current $\eta$ is the measure
  $\eta(B)$ of the base-ball $B$ of $\partial^2F$.
\end{defi}
The length of a current $\eta$ is explicitly given by
$$L(\eta)=\sum_{x\in A}\eta(\cyl{[1,x]}).$$
 
\begin{defi}[Frequency measures]
   The space of {\em frequency measures} is the set of finite-mass
  $T$-invariant 
  non-negative Borel measures on $\partial F$ (where $T$-invariant means that
  $T_*\mu=\mu$.) The total mass
  of a measure $\mu$ will be denoted by $||\mu||$. 
\end{defi}
The unitary ball 
  of the frequency measures, that is to say, the set of probability
  $T$-invariant measures on $\partial F$, is sometimes called the
  {\em frequency space} of $F$ in the literature (see~\cite{kap05}.)

 We will use the letter $\eta$ primarily to denote 
  a current, and the letter $\mu$ to denote a frequency measure.
We refer the reader to Appendix~\ref{a} for some basic
facts about currents and measures.

If $\eta$ is a geodesic current, and $x,y\in F$, by $F$-invariance, the
value $\eta(\cyl{[x,y]})$ depends only on the label $x^{-1}y\in
F$. The $F$-invariance of currents plays the role of $T$-invariance for
frequency measures. 
With this in mind, we can construct an isomorphism between the
space of geodesic currents and the 
space of frequency measures as follows: 
$$\eta \leftrightarrow\mu\qquad\textrm{ if and only if }\qquad 
\eta(\cyl{[x,y]})=\mu(\cyl{x^{-1}y}).$$
More precisely, one can prove (see also~\cite{kap05,kapre05}:)
\begin{teo}
The map $\alpha$ from the space of frequency measures to the space of
geodesic currents defined by
$$\alpha(\mu)(\cyl{[x,y]})=\mu(\cyl{x^{-1}y})$$
for all $x,y\in F$,
is a homeomorphism with respect to the weak-* topologies.
Moreover, under this correspondence, the total mass corresponds to length, that
is
$$L(\alpha(\mu))=||\mu||.$$
\end{teo}
\proof We only sketch the proof. The fact that $\alpha$ is
well-defined and bijective can be easily proved using $F$- and
$T$-invariance. The weak-$*$ continuity follows from
Proposition~\ref{p2.3_2.3}, while the last claim is a straightforward
computation.\qed

\

The identification between currents and frequency measures induces an
 action of Aut$(F)$ on the frequency measures given by 
$$\Phi\mu=\alpha^{-1}\circ(\Phi\times\Phi)_*\circ\alpha\mu.$$ 
Note that the action
on frequency measures is not just the push-forward via $\Phi$
because the push-forward does not commute with $T$.

The fact that length of a current corresponds to the total mass of a frequency 
measure will be the first ingredient of the proof of the compactness
result: bounded 
length $\to$ bounded norm $\to$ weak compactness.

\

Next, we briefly discuss relations between currents and cyclically
 reduced words, referring the reader to~\cite{kapre04}
for more details.

There is a natural embedding of the space of cyclically reduced words
in the space of geodesic currents (or frequency measures). Namely, if
$w$ is a cyclically reduced word, we denote by $w^{+\infty}$ the ray
$www\cdots$,  by $w^{-\infty}$ the ray $w^{-1}w^{-1}w^{-1}\cdots$, and 
by $\gamma_w$ the geodesic joining $w^{-\infty}$ and $w^{+\infty}$,
that is $\gamma_w=(w^{-\infty},w^{+\infty})\in\partial^2F$.
Then one can associate at each word $w$ the current
$$\eta_w=\sum_{v\in [w]}\delta_{\gamma_v}$$
where $[w]$ is the conjugacy class of $w$ in $F$ and
$\delta_{\gamma_v}$ 
denotes the Dirac measure concentrated on $\gamma_v$.
In the literature, such currents are often referred to as {\em rational
  currents}. Note that if $w$ is not a proper power, then
$||w||=L(\eta_w)=||\alpha^{-1}(\eta_w)||$.

\begin{defi}[Uniform current and uniform measure]\label{d3.1_3.5}
The {\em uniform current} $\eta_A$ and the 
{\em uniform frequency measure} $\mu_A$ are defined  
as follows. For all $v\in F$ we set
$$\mu_A(\cyl v)=\frac{1}{2k(2k-1)^{|v|-1}}
\qquad\textrm{ and }\qquad
\eta_A=\alpha(\mu_A).
$$
\end{defi}
Note that $L(\eta_A)=1$ and $||\mu_A||=1$.

\begin{remark}
  The uniform current is not the product
$\mu_A\times\mu_A$ on $(\partial F)^2$. Indeed, $\eta_A$ is a measure
on $\partial^2F\neq(\partial F)^2$, and $F$-invariance implies that
neighbourhoods of the diagonal have infinite measure. 
\end{remark}

Nevertheless, the current $\eta_A$ is not so different from
$\mu_A\times\mu_A$. Indeed, we can {\em disintegrate} $\eta_A$ with
measures that are in the same density class as $\mu_A$.
This means that if we cut a slice $S_x$ of $\partial^2F$ at the point $x$,
namely $S_x=\{x\}\times\{\partial F\setminus x\}$, then there exists a
continuous function $\f$ on $\{\partial F\setminus x\}$ such that the measure
$\mu_x$ induced on $S_x$ by $\eta_A$ is $\f\cdot\mu_A$.
A precise version of this fact is 
proved in Lemma~\ref{l2.2_2.4}.

\begin{defi}[Length of automorphisms]\label{defL}
For any automorphism $\Phi$ of $F$ we define the length of $\Phi$ as
the length of the image of the uniform current, that is
$$L(\Phi)=L(\Phi\eta_A)=\eta_A(\Phi^{-1}(B)).$$
Because of  $F$-invariance of currents, $L(\Phi)$ depends only on the
class $[\Phi]\in$Out$(F)$. We set $L([\Phi])=L(\Phi)$. 
\end{defi}

Intuitively speaking, for $\Phi\in\rm{Aut}(F)$, one can think of the
length of $\Phi$ as the limit, as $n$ tends to infinity, of
$||\Phi(w_n)||/n$ where $w_n$ is a ``random'' cyclically reduced word
in $F$ of length $n$.

\section{Proofs of the main results}\label{s4}

\proof[Sketch of the proof of Theorem~\ref{t2}] Let $\{\Phi_n\}$ 
be a sequence of
automorphisms of bounded length. The strategy of the proof can be
summarised as follows.

\vskip1ex
\step
The bounded length hypothesis, together with compactness of frequency
measures implies that the currents $\Phi_n\eta_A$ have a limit
$\eta_\infty$ (Lemma~\ref{l_step1}.) 

\vskip1ex\step The core of the proof.
We study of the action of $\Phi_n$ on $\partial^2F$ and on
$\partial F$. The main idea is that unbounded
lengths of $\Phi_n$ correspond to the fact that the maps $\Phi_n$
accumulate all  the boundary on some points (the {\em attractors}.) 
The key point is that, on the one hand, the bounded length hypothesis
excludes the presence of attractors, while on the other hand,
the existence of an element
of $f\in F$ such that $||\Phi_n(f)||$ is unbounded implies the presence of
attractors (Lemma~\ref{l4} and Lemma~\ref{l3}.) Therefore, the maps
 $\Phi_n$ keep bounded the cyclically reduced 
 length of all elements of $F$. 

\vskip1ex\step If $\Phi_n$ keep bounded the cyclically reduced 
 length of all elements of $F$, then
 $\Phi_n$ has a subsequence that converges, {\em i.e.} is eventually
 constant (Lemma~\ref{l_f}.)
 
\

We now proceed to work out the details of the proof.

\begin{lemma}\label{l_step1}
  Up to passing to a subsequence, the currents $\Phi_n\eta_A$ have a limit
  $\eta_\infty$ which is a geodesic current. 
\end{lemma}
\proof Since the lengths $L(\Phi_n\eta_A)$ are bounded, the total mass
of the corresponding frequency measure $\Phi_n\mu_A$ is bounded. 
The set of non-negative measures with bounded mass on a compact space 
is weakly-* compact. Since $\partial F$ is compact, up to passing to
subsequences, $\Phi_n\mu_A$ has a limit $\mu_\infty$. Such a limit is
$T$-invariant because the push-forward via a continuous map 
is weak-* continuous. Thus,
the limit is a frequency measure, which therefore corresponds to a geodesic
current $\eta_\infty$. By continuity of the correspondence $\alpha$ between
frequency measures and geodesic currents, it follows that
$\Phi_n\eta_A\to \eta_\infty$.\qed

\begin{remark} Although it is not relevant for the sequel, we note
  that the limit $\eta_\infty$ is actually a non-zero current because
  $L(\Phi)\geq1$ for any $\Phi$.
\end{remark}

Next, we proceed to study the attractors. As the next lemma shows,
attractors correspond to singularities of the limit current
$\eta_\infty$. More precisely,
a Borel measure $\sigma_1$ is said to be {\em absolutely continuous}
w.r.t. the Borel measure  $\sigma_2$ if for any 
Borel set $C$, $\sigma_2(C)=0$ implies $\sigma_1(C)=0$. We say
that a current $\eta$ {\em has a part concentrated on a set} $C$ if 
$\eta_A(C)=0$ and $\eta(C)>0$. Similarly, we say that
a frequency measure $\mu$ {\em has a part concentrated on a set} $C$ if 
$\mu_A(C)=0$ and $\mu(C)>0$.

\begin{lemma}\label{l2.3_3.6}
Let $\Psi_n$ be a sequence of automorphisms such that the currents
$\Psi_n\eta_A$ have a limit current $\eta$. 
Let $p,q\in\partial F$ be two distinct points.
Suppose that there exist cylinders
$P_n=\cyl{p_n},Q_n=\cyl{q_n}\subset\partial F$ such that $p_n\to p$
and $q_n\to q$, and such that
there is a positive constant $c$ for which
$\eta_A(\Psi_n^{-1}(P_n\times Q_n))>c$. Then,
the current $\eta$ has a part concentrated on
$(p,q)$.
\end{lemma}
\proof
Any cylinder $C\subset\partial^2 F$ containing $(p,q)$, 
contains $P_n\times Q_n$ for all sufficiently large $n$. 
Therefore, by definition of
push-forward and by hypothesis, we get 
(for all sufficiently large $n$:)
$$\Psi_n\eta_A(C)>\Psi_n\eta_A(P_n\times
Q_n)=\eta_A(\Psi_n^{-1}(P_n\times Q_n))>c$$
By Proposition~\ref{p2.3_2.3} it follows that the limit current
satisfies $\eta(C)\geq c$ for any cylinder $C$ containing $(p,q)$.
This implies that $\eta((p,q))\geq c$ while $\eta_A((p,q))=0$, that is,
$\eta$ has a part concentrated on $(p,q)$.\qed

\

The following is a standard fact, which says that the only currents that
can have a part concentrated on a geodesic are essentially 
 the rational currents (recall that $\partial F$ and $\partial^2 F$
 are identified with the set of geodesic rays and of
 bi-infinite geodesics in the Cayley graph respectively.)

\begin{lemma}\label{l1} Any frequency measure (and hence
  the limit $\mu_\infty$) cannot have a part concentrated on a non-periodic
  ray. Therefore any current (and hence the limit 
$\eta_\infty$) cannot have part concentrated on a
  non-periodic geodesic.
\end{lemma}
\proof Any frequency  measure $\mu$ has finite mass and is
$T$-invariant. Hence if it has a part
concentrated on a point $x$, it has a part concentrated on each point of
the $T$-orbit of $x$, with the same weight. 
It follows that such an orbit must be finite, forcing $x$ to be
periodic.
\qed

\

We prove now a small lemma that will be used in Lemma~\ref{l4}. 
We isolated this result  
from the proof of Lemma~\ref{l4} because it is the only point in
which we crucially use the properties of the uniform current
(this will be further discussed in Section~\ref{s5}.) 

\begin{lemma}\label{l4344} For any $f\in F$, there exists  a positive
  constant $c$, depending only on $f$, such that whenever 
two disjoint Borel-subsets $E,S$ of $\partial F$ satisfy
$$f(\partial F\setminus E)\subset S\qquad f^{-1}(\partial F\setminus
S)\subset E,$$
then we have 
$$\eta_A(E\times S)\geq c(1-\mu_A(E))(1-\mu_A(S)).$$
In particular, $\eta_A(E\times S)=0$ if and only if
either $E$ or $S$ has full-measure with respect to $\mu_A$.
\end{lemma}
We notice that the hypothesis that $E$ and $S$ are disjoint is
redundant. We formulated the lemma in that way for reasons of
compatibility with Theorem~\ref{ttt}.
 
\proof[Proof of Lemma~\ref{l4344}]
The hypothesis that $E$ and $S$ are disjoint will not be used in this proof.
By hypothesis we have $\mu_A(S)\geq\mu_A(f(\partial F\setminus E))$ 
and $\mu_A(E)\geq\mu_A(f^{-1}(\partial F\setminus S))$, 
and by Lemma~\ref{l2.2_2.6}
$$\mu_A(E)
\geq\frac{1-\mu_A(S)}{(2k-1)^{|f|}}$$
and
$$\mu_A(S)
\geq\frac{1-\mu_A(E)}{(2k-1)^{|f|}}$$
from which we get
$$
\mu_A(E)\mu_A(S)\geq
\frac{(1-\mu_A(E))(1-\mu_A(S))}{(2k-1)^{2|f|}}
.
$$

By Lemma~\ref{l2.2_2.5} 
$$
\eta_A(E\times S)\geq\mu_A(E)\mu_A(S)
\geq
\frac{(1-\mu_A(E))(1-\mu_A(S))}
{(2k-1)^{2|f|}}
$$
and setting $c=\frac{1}{(2k-1)^{2|f|}}$ completes the proof.
\qed

\

Now the aim is to prove that the maps $\Phi_n$ keep bounded the
lengths of all elements of $F$, so that we can apply Lemma~\ref{l_f}.
We do it in the following two lemmata. Namely, 
in Lemma~\ref{l4} we show that if this is not the case, then there
are no attractors in $\partial^2F$ and at most a unique attractor in
$\partial F$. Lemma~\ref{l3} will show that, up to conjugations, we can
avoid the presence of a unique attractor in $\partial F$.

\begin{lemma}\label{l4} Suppose that $\Phi_n\eta_A$ has bounded
  length (uniformly on $n$.)
  Suppose that there exists an element $f\in F$ such that the
  cyclically reduced length $\Phi_n(f)$ goes to $\infty$.
Then, after possibly passing to a subsequence, 
$\Phi_n$ (as maps of $\partial F$) pointwise converge almost everywhere
  to a  constant. That is to say, 
  there exists $y\in\partial F$ such that for $\mu_A$-almost all
  $x\in\partial F$, $\Phi_n(x)\to y.$ 
\end{lemma}
\proof
By Lemma~\ref{l_step1}, without loss of generality we can
suppose that 
$\Phi_n\eta_A$ has a limit $\eta_\infty$.

Recall that we consider $\Phi_n(f)$ and $\Phi_n(f^{-1})$ as freely
reduced words. 
Let $v_n$ be the maximal
initial segment shared by $\Phi_n(f)$ and $\Phi_n(f^{-1})$, and let
$\Psi_n$ be the map $x\mapsto v_n^{-1}\Phi_n(x) v_n$. Note that
$\Psi_n\eta_A=\Phi_n\eta_A$ and that
$\Psi_n(f)$ is cyclically reduced.
Up to passing to a subsequence,
$\Psi_n(f)$ and $\Psi_n(f^{-1})$ have limits, which we denote by $r_+$ and
$r_-$, in $\partial F$. Since $\Psi(f)$ is cyclically reduced,
$r_+\neq r_-$. Note that 
this also implies that $r_+$ and $r_-$ have no common initial segment,
that is, the geodesic $(r_-,r_+)$ passes through $1$, the base-point
of the Cayley graph. 

We now show that $\Psi_n$, as maps of $\partial F$, converges $\mu_A$-almost
everywhere either to $r_-$ or to $r_+$.

Next, cut $\Psi_n(f)$ into two segments of equal
length. More
precisely, we set 
 $$\Psi_n(f)=s_ne_n^{-1}$$
where the starting segment $s_n$ and the ending one $e_n$ both have
length $|\Psi_n(f)|/2$ (approximated to the nearest integers.) 
We have $s_n\to r_+$ and $e_n\to r_-$. In particular, for large enough
$n$, $\cyl{e_n}\cap\cyl{s_n}=\emptyset$, which implies
$\cyl{[e_n,s_n]}=\cyl{e_n}\times\cyl{s_n}$. For large $n$, 
let $C_n$ be such a cylinder:
$$C_n=\cyl{[e_n,s_n]}=\cyl{e_n}\times\cyl{s_n}$$
and set $E_n=\Psi_n^{-1}(\cyl{e_n})$ and $S_n=\Psi_n^{-1}(\cyl{s_n})$.
Note that $E_n\cap S_n=\emptyset$.

For all $x\in\partial F$, either
$\Psi_n(x) \in \cyl{e_n}$ or $\Psi_n(fx) \in \cyl{s_n}$, so
either $x \in E_n$ or 
$fx \in S_n$,
whence 
$$f(\partial F\setminus E_n)\subset S_n.$$
Similarly, either $\Psi_n(x) \in \cyl{s_n}$ or
$\Psi_n(f^{-1}x) \in \cyl{e_n}$ and
$$f^{-1}(\partial F\setminus S_n)\subset E_n.$$ 
Thus, by Lemma~\ref{l4344}
$$
\eta_A(E_n\times S_n)\geq
\frac{(1-\mu_A(E_n))(1-\mu_A(E_n))}{(2k-1)^{2|f|}}
.
$$

By definition of push-forward
$$
\Psi_n\eta_A(C_n)=\eta_A(\Psi_n^{-1}(C_n))=\eta_A(E_n\times S_n),
$$
and putting together these (in)equalities, we get
$$
\Psi_n\eta_A(C_n)
\geq
\frac{\big[1-\mu_A\big(\Psi_n^{-1}(\cyl{e_n})\big)\big]\big[1-\mu_A\big(\Psi_n^{-1}(\cyl{s_n})\big)\big]} 
{(2k-1)^{2|f|}}
. 
$$

If $\Psi_n\eta_A(C_n)\to 0$ then either 
$\mu_A(\Psi_n^{-1}(\cyl{e_n}))$ or
$\mu_A(\Psi_n^{-1}(\cyl{s_n}))$ converges to $1$ and therefore,
up to passing to subsequences,
 $\Psi_n$ converges almost
everywhere either to $r_-$ or to $r_+$, and we are done.

We now show that the bounded length hypothesis excludes the
possibility that
$\Psi_n\eta_A(C_n)$ stays bounded away from zero. Indeed, suppose
that there exists a constant $c$ such that
$\Psi_n\eta_A(C_n)>c$, uniformly on $n$. Then,  by Lemma~\ref{l2.3_3.6},
$\eta_\infty$ has a part concentrated on 
the geodesic $(r_-,r_+)$, which is therefore periodic 
by Lemma~\ref{l1};  
let $w$ be its period. We must have 
\begin{equation}\label{e_num}
r_-=w^{-\infty} \qquad \textrm{ and } \qquad r_+=w^{+\infty}.
\end{equation}
We may assume that the element $f$ is not a proper power. Since
$\Psi_n$ is an automorphism of $F$, it follows that $\Psi_n(f)$ is not
a proper power either.
Therefore, by~$(\ref{e_num})$, for all large enough $n$, 
we can write $\Psi_n(f)$ as
$$\Psi_n(f)=w^{i(n)}u_nw^{j(n)}$$
 where $u_n$ neither
starts nor ends with $w$, and the exponents $i(n)$ and $j(n)$ are
non-negative and go to infinity as $n$ does. Without loss of
 generality, we can suppose $i(n)\leq j(n)$, so that $s_n$ starts
with  $w^{i(n)}$.

For any $0\leq h\leq i(n)$, 
let $C_n^h$ be the cylinder
$$C_n^h=\cyl{[w^{-h}e_n,w^{-h}s_n]}.$$ 
Note that the $C_n^h$'s are pairwise disjoint, because $u_n$ neither
starts nor ends with $w$. Moreover,   
the condition $0\leq h\leq i(n)\leq j(n)$ implies that
the geodesic segment from $w^{-h}e_n$ and $w^{-h}s_n$ passes
through $1$, thus 
$C_n^h\subset B$ for all $0\leq h\leq i(n)$.
By $F$-invariance of currents, we have
$$\Psi_n\eta_A(C_n^h)=\Psi_n\eta_A(C_n)>c$$
uniformly on $n$. It follows that
$$L(\Phi_n)=L(\Psi_n)=\Psi_n\eta_A(B)>ci(n)$$ 
which goes to infinity as $n$ does, contradicting the bounded length
hypothesis. 

Thus we have proved that, after passing to a subsequence, the maps
$\Psi_n:x\mapsto v_n^{-1}\Phi_n(x)v_n$ $\mu_A$-almost everywhere converge to a
map which is constant (either to $r_-$ or to $r_+$). 
Up to possibly passing
to a subsequence, $v_n$ converges to a limit $v_\infty$, which is
either an element of $F$ or of $\partial F$.
Since the elements $v_n$ were the maximal
initial segments shared by $\Phi_n(f)$ and $\Phi_n(f^{-1})$, the words
$v_n\Psi_n(f)$ and $v_n\Psi_n(f^{-1})$ are freely reduced.
It follows that the maps $\Phi_n$ converge almost everywhere, 
up to passing to the same
subsequence, 
to a constant -- which is either $v_\infty$ (if $v_\infty$ is an
element of $\partial F$) or $v_\infty r_-$ or $v_\infty r_+$ (if
$v_\infty\in F$.) \qed

\begin{remark}
  The proof of Lemma~\ref{l4} can be adapted to prove the following
  more general fact. 
  If we replace the hypothesis ``$\Phi_n\eta_A$ has bounded length
  $\eta_\infty$'' 
  with ``$\frac{\Phi_n}{L(\Phi_n)}\eta_A$ has a limit $\eta_\infty$''
   -- which is always true up to passing to a subsequence -- then, we get that
  $\eta_\infty$ does not have a part concentrated on a geodesic. Indeed, if
  $\eta_\infty$ has a part concentrated on a geodesic $\gamma$, then
  there exists a positive constant $c$ such that
  for any cylinder $C$ containing $\gamma$ we have 
  $\eta_A((\Phi_n\times\Phi_n)^{-1}(C))\sim cL(\Phi_n)$.
  As in the argument above, we must have  
  $\gamma=(w^{-\infty},w^{+\infty})$ for some $w\in
  F$, and conjugating $\Phi_n$ by a suitable power of $w$, we reach a
  contradiction. Indeed, if $X$ denotes
  the set  $(\Phi_n\times\Phi_n)^{-1}(C)$, then 
  $(\Phi_n\times\Phi_n)^{-1}(wC)=\Phi_n^{-1}(w) X$ which is contained
  in the set $(\Phi_n\times\Phi_n)^{-1}(B)$, whose $\eta_A$-measure is
  $L(\Phi_n)$ by definition. If a geodesic belongs to
  $X\cap\Phi_n^{-1}(w)X$, then it passes through $1$, whence
  $\eta_A(X\cap\Phi_n^{-1}(w)X)\leq 1$. Since $\eta_A(X)\sim
  cL(\Phi_n)$, we can
  conjugate by $w$ approximately at most $1/c$ times, while 
 if $\gamma=(w^{-\infty},w^{+\infty})$, then we can do that infinitely many
  times.  
\end{remark}

After Lemma~\ref{l4}, it remains to deal with the case where
$\Phi_n$ converges almost 
everywhere to a constant.
What is the behavior of such a sequence? An example can be constructed by
taking a fixed $\Phi$ and conjugating with elements $v_n$ whose
length goes to infinity. The next lemma shows
that more or less this is the only possibility.

\begin{lemma}\label{l3}
Let $\Phi_n$ be a sequence of automorphisms of $F$. 
 Then, there exists $v_n\in F$ such that, up possibly
  passing to a subsequence, the maps $x\mapsto v_n^{-1}\Phi_n(x)v_n$
  have no subsequence converging to a constant $\mu_A$-almost everywhere.
\end{lemma}
\proof The rough idea is that, via conjugations, we can force the
``barycentre of $\Phi_n$'' to stay in a fixed compact.

For any freely reduced word $w$ of length $M$, define $B_n(w)$ as the
set of rays $x$  
such that $\Phi_n(x)$ starts by $w$, namely
$$B_n(w)=\{x\in\partial F:\Phi_n(x)\in\cyl w\}=\Phi_n^{-1}(\cyl w).$$

Obviously $B_n(1)=\partial F$. Moreover,
 for each $n$ we have:
 \begin{equation}\label{e2}
\lim_{M\to\infty}\sup_{|w|=M}\mu_A(B_n(w))=0.   
 \end{equation}
Indeed, otherwise for all $M$ there exists $w_M\in F$ of length $M$
 such that $\mu_A(B_n(w_M))>c>0$.
Up to subsequences, $w_M$ converges to a ray $R$, and
$\Phi_n(B_n(w_M))=\cyl{w_M}\to R$, contradicting the fact that 
$\Phi_n$ is a homeomorphism of $\partial F$ (in this argument $n$ is fixed.)

Now, let $v_n$ be a freely reduced word of maximal length such that
$\mu_A(B_n(v_n))\geq \frac{1}{2}$.
Let $\wt\Phi_n$ be the map $x\mapsto v_n^{-1}\Phi(x)v_n$ and let
$$\wt B_n(w)=\{x\in \partial F:\wt\Phi_n(w)\in\cyl w\}.$$
Let $l_n\in A$ be the last letter of $v_n$. Since
$\mu_A(B_n(v_n))\geq\frac{1}{2}$ we get $\mu_A(\wt
B_n(l_n^{-1}))\leq\frac{1}{2}$. On the other hand, for any $a\in A$, different 
from $l_n^{-1}$, maximality of the length of $v_n$ implies 
$$\mu_A(\wt B_n(a))\leq\frac{1}{2}.$$
Hence, such an inequality holds for all $a\in A$.
It follows that the sequence $\wt\Phi_n$ cannot have any subsequence
converging to a constant almost everywhere.\qed

\

Since conjugations do not affect the length of cyclically reduced
words, Lemma~\ref{l4} and Lemma~\ref{l3} can be summarised as
follows (recall that for $f\in F$, $|f|$ denotes its length, while
$||f||$ denotes the  length of the 
cyclically reduced word obtained from $f$.)

\begin{cor}\label{c_a_4.9}
  Let $\Phi_n$ be a sequence of automorphisms. If there is
  $M$ such that  
  $L(\Phi_n)<M$, then for each $f\in F$ there exists $M(f)$ such that
  $||\Phi_n(f)||<M(f)$. 
\end{cor}

As the experts know, Corollary~\ref{c_a_4.9} implies Theorem~\ref{t2}. \
We include the proof of the following 
Lemma~\ref{l_f} for completeness.  

\begin{lemma}\label{l_f}
  Let $\{\Phi_n\}$ be a sequence of automorphisms such that 
  for each $f\in F$ there is an $M(f)$ such that
  $||\Phi_n(f)||<M(f)$. Then, there exist elements $v_n\in F$ 
  such that a subsequence of $\{v_n^{-1}\Phi_nv_n\}$ converges to an
  automorphism ({\em i.e.} $\{v_n^{-1}\Phi_nv_n\}$ has a constant
  subsequence.) 
\end{lemma}

\proof 
By a diagonal argument, up to passing to a subsequence, the
maps $\Phi_n$, as maps from $F$ to itself, pointwise converge to a
map $\Phi_\infty$ (up to conjugations.) In particular, there exists a
map $\Phi_\infty:A\to F$ and maps $w_n:A\to 
F$ such that $\Phi_\infty(f)$ is cyclically reduced and, up to
passing to a subsequence, for all $f\in A$ we have 
$$\Phi_n(f)=w_n(f)\Phi_\infty(f)w_n(f)^{-1}.$$

Choose an element $a\in A$. Up to conjugations we can suppose that
$w_n(a)=1$, that is, $\Phi_n$ really converges as an automorphism on
the subgroup 
generated by $a$. Let $G\subset A$ be a maximal set of generators $g$ such that
$|\Phi_n(g)|$ stays bounded. If $G=A$ we are done, because, up to
subsequences, $\Phi_n$ converges on $A$, whence on $F$. Otherwise, there
exists $f\in A$ such that the length of $w_n(f)$ goes to infinity. 
Since
$$\Phi_n(af)=\Phi_\infty(a)\Phi_n(f)=\Phi_\infty(a)w_n(f)\Phi_\infty(f)w_n(f)^{-1}$$
has bounded cyclically
reduced length, and since $\Phi_\infty(a)$ has finite length,
we get that, for large enough $n$,
 $w_n(f)$ must start either with $\Phi_\infty(a)$ or with
$\Phi_\infty(a)^{-1}$. Iterating this argument we get that
$w_n(f)$ is the product of a power of $\Phi_\infty(a)$ and a word of
bounded length.  Thus, up to subsequences, we get  
\begin{equation}\label{eq2}w_n(f)=\Phi_\infty(a)^{m}u\end{equation} 
for some $m\in\M Z$  with $|m|\to\infty$ as $n\to\infty$, 
and $u$ a finite word (which depends on $f$.)

It follows that, up to conjugating $\Phi_n$ by
$\Phi_\infty(a)^m$, we can suppose that $G$ has at least two elements
$a,b$ and that $\Phi_n$ is eventually constant on the subgroup
generated by $a$ and $b$.
If $G\neq A$, let $f$ be as above.
As in $(\ref{eq2})$, we get
$$w_n(f)=\Phi_\infty(a)^{m}u$$
$$w_n(f)=\Phi_n(b)^{l}v$$
for some exponents $m,l$ such that $|m|,|l|$ go to
infinity as $n$ does, and fixed words $u,v$ (depending on $f,a,b$.)
Therefore, as $n$ goes to infinity, we get that the unoriented geodesics
$(\Phi_\infty(a)^{-\infty},\Phi_\infty(a)^{+\infty})$ and 
$(\Phi_n(b)^{-\infty},\Phi_n(b)^{+\infty})$ coincide.
 
  This implies that $\Phi_n(b)$ is cyclically
reduced. In particular, we get  
$\Phi_n(b)=\Phi_\infty(b)$, and therefore $\Phi_\infty$ is an
automorphism on the group generated by $a$ and $b$. Moreover, the above
inequalities imply that 
$$\Phi_\infty(b)^{|\Phi_\infty(a)|}=\Phi_\infty(a^{\pm1})^{|\Phi_\infty(b)|}$$
whence
$$b^{|\Phi_\infty(a)|}=a^{\pm1|\Phi_\infty(b)|}$$
which is impossible because $F$ is free. Thus $G=A$, and hence there exists a
subsequence of $\{\Phi_n\}$ which converges.
\qed

\vskip\baselineskip
By Lemma~\ref{l3}, up to conjugations, the sequence $\Phi_n$ does not
sub-converge almost everywhere to the same point; by Lemma~\ref{l4}
we can apply Lemma~\ref{l_f}, and the proof of Theorem~\ref{t2} is
complete.\qed

\proof[Proof of Corollary~\ref{c2.1_2.2}] We have to prove that the
spectrum of the length function is discrete.
Suppose not, and take a sequence $\Phi_n$ of automorphisms
such that $L(\Phi_n)$ has a limit $\lambda$, with
$L(\Phi_n)\neq\lambda$  for all $n$.
By Theorem~\ref{t2} there exist elements $v_n$ and a subsequence $n_i$
such that the maps $\Psi_{n_i}:x\mapsto v_{n_i}\Phi_{n_i}(x)v_{n_i}^{-1}$
converge to an automorphism $\Psi$. Thus, the sequence $\Psi_{n_i}$
is eventually constant, and therefore the sequence of lengths
$L(\Psi_{n_i})$ is also eventually constant. But
$L(\Psi_{n_i})=L(\Phi_{n_i})$ is therefore
eventually equal to $\lambda$, a contradiction.\qed

\proof[Proof of Theorem~\ref{t_2.4_1.3}]
This immediately follows from Corollary~\ref{c2.1_2.2} 
and \cite[Proposition~5.2]{kapre05}. Indeed, I. Kapovich proved that
for any non-simple automorphism $\Phi$ there exists a Whitehead
automorphism $\tau$ such that 
$$1\leq L(\tau\Phi)<L(\Phi)$$
and the claim follows by an inductive argument on the length.\qed

\section{Generalisations}\label{s5}

In this section we give a (partial) answer to the question: {\em For
  which currents does Theorem~\ref{t2} hold?}. 
The idea is that length-compactness is true
(for any action on metric trees and) for any current for which
Lemma~\ref{l4344} holds.

\begin{defi}[$\eta$-length of automorphisms]
Let $\eta$ be a geodesic current and let
$\Phi\in\textrm{Aut}(F)$. We define the $\eta$-length of $\Phi$ as
$$L_\eta(\Phi)=L(\Phi\eta).$$
\end{defi}

The proof of Theorem~\ref{t2} can be followed step by step in
this new setting, obtaining:

\begin{teo}\label{t1}
  Let $\eta$ be a geodesic current and let $\mu$ be its
  corresponding frequency-measure. Suppose there exists a
  constant $c>0$ and that for each $f\in F$ there is $b(f)>0$ such
  that for any disjoint  Borel sets $E,S\subset\partial F$ 
$$\mu(fE)\geq b(f)\mu(E)\qquad\textrm{ and } \qquad \eta(E\times
S)\geq c\mu(E)\mu(S).$$
Then Theorem~\ref{t2} holds for $\eta$. That is to say, any sequence of
automorphisms $\Phi_n\in\textrm{Aut}(F)$ with bounded $\eta$-length
has, after possibly conjugating, a convergent subsequence ({\em i.e.}
a constant subsequence.)
\end{teo}
\proof
The hypotheses on $\eta$ guarantee that Lemma~\ref{l4344} holds for
$\eta$. Moreover, our assumptions imply that $\eta$ is not
concentrated on a single geodesic. This implies that Lemma~\ref{l3}
can be rewritten, with the difference that equation~$(\ref{e2})$ of
page~\pageref{e2} becomes (notation as in Lemma~\ref{l3})
$$\lim_{M\to\infty}\sup_{|w|=M}\mu(B_n(w))<C<1.$$
for a certain constant $C$, so that we have to consider a word $v_n$ of
maximal length such that $\mu(B_n(v_n))\geq C$. 

The proof now is exactly as in Theorem~\ref{t2} because 
there are no other places in the proof of Theorem~\ref{t2} where
the specific properties of the uniform current were used.\qed

\

We now give a criterion for a current to satisfy the
hypotheses of Theorem~\ref{t1}, formulated in terms of the corresponding
frequency measure.

\begin{teo}\label{ttt} Let $\eta$ be a geodesic current and $\mu$ be its
  corresponding frequency measure.
  Suppose that for each $a\in A$ there exist two strictly positive constants
  $C_1(a)$ and $C_2(a)$ such that for any $E\subset \partial
  F\setminus\cyl{a^{-1}}$ we have
$$C_1(a)\mu(E)\leq\mu(aE)\leq C_2(a)\mu(E).$$
For all freely reduced word $w=a_0\dots a_k$,  set $C_i(w)=C_i(a_0)\cdots
C_i(a_k), i=1,2$.   
  
If there is a constant $M$ such that
$$\inf_{w\in F}\frac{C_1(w^{-1})}{C_2(w)}\geq\frac{1}{M}\qquad\sup_{w\in F}
C_2(w)\leq M$$ 
then the hypothesis of Theorem~\ref{t1} is satisfied. Namely,
there exists a
  constant $c>0$ and for each $f\in F$ there is $b(f)>0$ such
  that for any $E,S\subset\partial F$ 
$$\mu(fE)\geq b(f)\mu(E)\qquad\textrm{ and } \qquad \eta(E\times
S)\geq c\mu(E)\mu(S).$$
In particular, length compactness holds for $\eta$.
\end{teo}
\proof By Proposition~\ref{p2.3_2.3}, it suffices to prove the claims
when $E$ and $S$ are cylinders, and, regarding the first claim, it
suffices to prove it for generators.
Let $E\subset \partial F$ be a
cylinder and $a\in A$. Set $E_0=E\cap\cyl{a^{-1}}$ and
$E_1=E\setminus E_0$, then
$$\mu_A(aE_0)\geq\frac{1}{C_2(a^{-1})}\mu_A(E_0)
\qquad\mu_A(aE_1)\geq C_1(a)\mu_A(E_1)$$ 
and the first claim follows setting
$b(a)=\inf\{C_1(A),\frac{1}{C_2(a^{-1})}\}$. 

We now prove the second claim. Let $E,S\subset\partial F$ be two
disjoint cylinders. 
Let 
$$E=\cyl{v_0v}\qquad S=\cyl{v_0w}$$
with  $v_0v,v_0w,v^{-1}w$ freely reduced (that is, $v_0$ is the
maximal initial segment shared by $E$ and $S$.) We set
$$E'=\cyl v=v_0^{-1}E\qquad S'=\cyl w=v_0^{-1}S$$
$$E''=v^{-1}E'\qquad S''=v^{-1}S'=\cyl{v^{-1}w}$$

By $F$-invariance we have
$$\eta(E\times S)=\eta(E'\times S')=\eta(E''\times S'')=\eta(\cyl{[v^{-1}w]})=\mu(S'')$$

Then by induction on the lengths of $v_0,v,w$ (and
using that $v_0v,v_0w,v^{-1}w$ are reduced) we get
\begin{eqnarray*}
  & &\mu(S'')\geq C\mu(S'')\mu(E'')\geq
  CC_1(v^{-1})\mu(S')\frac{1}{C_2(v)}\mu(E')\\
&=&c\frac{C_1(v^{-1})}{C_2(v)}\mu(S')\mu(E')\geq
C\frac{C_1(v^{-1})}{C_2(v)}
\frac{\mu(S)}{C_2(v_0)}\frac{\mu(E)}{C_2(v_0)}
\end{eqnarray*}
where $C$ is a suitable constant, and the claim follows by setting
$c=C/M^3$.\qed

\appendix
\section{}\label{a}
Throughout the paper, we used some standard results about currents and
measures. This section contains the proofs of some of these facts.

  \begin{prop}\label{p2.3_2.3}
    Let $m$ be a Borel measure on $\partial F$ or $\partial^2F$. 
    Then, $m$ is determined by its value on cylinders. 
    Moreover, if $\{m_i\}$ is
    a sequence of Borel measures, then  
    $m_i$ converges to $m$ if and only if for all cylinders $C$,
    $m_i(C)\to m(C)$.
    
  \end{prop}
\proof We restrict to the case where $m$ is 
a Borel measure on $\partial F$; an identical argument
gives the $\partial^2 F$-case.
So assume $m$ is a Borel measure on $\partial F$.
The characteristic function of any cylinder belongs to $\call
C_0(\partial F)$, and the space
$V$ generated by the characteristic functions of cylinders is dense in
$\call C_0(\partial F)$ (the topology  of $\call C_0(\partial F)$ is given
by uniform convergence.) The first claim follows.

In addition, 
this implies that if $m_i$ converges to $m$, then for any cylinder $C$, 
$m_i(C)\to m(C)$. 
On the other hand, suppose that $m_i(C)\to m(C)$ for all cylinders
$C$. Then for any $\chi\in V$, $\int\chi\,dm_i\to\int\chi\,dm$. 
Therefore, for any $\f\in\call C_0(\partial F)$, if $\{\chi_k\}\subset
V$ is a 
sequence converging to $\f$, we have:
$$|\int\f\,d(m_i-m)|\leq
|\int|\f-\chi_k|\,dm_i|+
|\int\chi_k\,d(m_i-m)|+
|\int|\chi_k-\f|\,dm|
$$
where the sum of the first and the last term is bounded by 
$||\f-\chi_k||(||m_i||+||m||)$, which goes to zero as $k\to\infty$, 
 uniformly on $i$.
The second term goes to zero for any $k$.
\qed

\begin{lemma}\label{l2.2_2.4} For any $(x,y)\in \partial^2 F$ let $L(x,y)$
  be the length of the maximal initial segment shared by $x$ and $y$.
  Then we have
$$\eta_A=2k(2k-1)^{2L(x,y)-1}\mu_A\times\mu_A.$$
\end{lemma}

\proof 
  Let $D,E\subset\partial F$ be two disjoint cylinders. 
  Since $D$ and $E$ are disjoint, there exist $v,w\in
F$ such that $D=\cyl v$, $E=\cyl w$, and such that $v$ is
not the initial segment of $w$ and vice versa. Let $L$ be the length of
the maximal initial segment shared by $v$ and $w$ (possibly $L=0$.)
Now, let $D'=\cyl{v'}\subset D$ and $E'=\cyl{w'}\subset E$ be two
cylinders. We set $|v'|=L+a$ and $|w'|=L+b$. We have
$D'\times E'=\cyl{[v',w']}$ and, by definition (Definition~\ref{d3.1_3.5})
$$\eta_A(\cyl{[v',w']})=\frac{1}{2k(2k-1)^{|(v')^{-1}w'|-1}}
=\frac{1}{2k(2k-1)^{a+b-1}}$$
which can be written as
$$\frac{2k(2k-1)^{2L-1}}{(2k(2k-1)^{L+a-1})(2k(2k-1)^{L+b-1})}
=2k(2k-1)^{2L-1}\mu_A(D')\mu_A(E)'$$
So by Proposition~\ref{p2.3_2.3}, the restriction of $\eta_A$ to
$D\times E$ is given by $2k(2k-1)^{2L-1}\mu_A\times\mu_A.$
Since for each $(x,y)\in D\times E$ we have $L(x,y)=L$, we get that
the restriction of $\eta_A$ to $D\times E$ is given by 
$$2k(2k-1)^{2L(x,y)-1}\mu_A\times\mu_A.$$
Since this holds for any $D,E$, 
the claim follows by Proposition~\ref{p2.3_2.3}.
\qed

\

An immediate corollary of Lemma~\ref{l2.2_2.4} is the following.
\begin{lemma}\label{l2.2_2.5}
Let $E,D\subset \partial F$ be two Borel subsets of $\mu_A$-positive
measure. Then $\eta_A(E\times D)\geq \mu_A(E)\mu_A(D)$.  
\end{lemma}
\proof
Just apply Fubini-Tonelli theorem, using Lemma~\ref{l2.2_2.4}, and the
fact that $L(x,y)\geq 0$ and that $\frac{2k}{2k-1}>1$.
\qed

\begin{lemma}\label{l2.2_2.6}
  Let $E$ be a Borel subset of $\partial F$. Then for all $f\in F$ 
$$\mu_A(fE)\geq\frac{\mu_A(E)}{(2k-1)^{|f|}}$$ 
In particular, if $E$ has $\mu_A$-positive
  measure, then $fE$ has $\mu_A$-positive measure.
  \end{lemma}
\proof It suffices to prove the claim for $f\in A$. 
Let $E_0=E\cap\cyl{f^{-1}}$ and $E_1=E\setminus E_0$. 
By definition of $\mu_A$ and Proposition~\ref{p2.3_2.3} we have
$$\mu_A(fE_0)=(2k-1)\mu_A(E_0)\qquad\mu_A(fE_1)=\frac{\mu_A(E_1)}{2k-1}$$
and the claim follows.\qed


\begin{thebibliography}{Kap05b}

\bibitem[AFP00]{afp:libro}
Luigi Ambrosio, Nicola Fusco, and Diego Pallara.
\newblock {\em Functions of bounded variation and free discontinuity problems}.
\newblock Oxford Mathematical Monographs. The Clarendon Press Oxford University
  Press, New York, 2000.

\bibitem[Bon86]{Bon86}
Francis Bonahon.
\newblock Bouts des vari\'et\'es hyperboliques de dimension {$3$}.
\newblock {\em Ann. of Math. (2)}, 124(1):71--158, 1986.

\bibitem[Coo87]{coo87}
Daryl Cooper.
\newblock Automorphisms of free groups have finitely generated fixed point
  sets.
\newblock {\em J. Algebra}, 111(2):453--456, 1987.

\bibitem[Fed69]{fed:libro}
Herbert Federer.
\newblock {\em Geometric measure theory}.
\newblock Die Grundlehren der mathematischen Wissenschaften, Band 153.
  Springer-Verlag New York Inc., New York, 1969.

\bibitem[Kap05]{kap05}
Ilya Kapovich.
\newblock The frequency space of a free group.
\newblock {\em Internat. J. Alg. Comput. (special {G}aeta {G}rigorchuk's 50's
  birthday issue)}, 15(5-6):939--969, 2005.

\bibitem[Kap06a]{kapre05}
Ilya Kapovich.
\newblock {C}lusters, currents and {W}hitehead's algorithm.
\newblock {\em Prerint, ar{X}iv:math.{GR}/0511478, to appear on 
Experimental Mathematics}, 2006.

\bibitem[Kap06b]{kapre04}
Ilya Kapovich.
\newblock {C}urrents on free groups.
\newblock {In \em Topological and asymptotic aspects of group theory}, 
volume 394 of {\em Contemp. Math.}, pages 149--176.
\newblock{Amer. Math. Soc., Providence, RI}, 2006.

\bibitem[KKS05]{kks06}
Vadim Kaimanovich, Ilya Kapovich, and Paul Shupp.
\newblock The subadditive ergodic theorem and generic stretching factors for
  free group automorphisms.
\newblock {\em Prerint, ar{X}iv:math.{GR}/0504105, to appear on Israel J.
  Math.}, 2005.

\bibitem[LS77]{lysc:libro}
Roger C. Lyndon and Paul E. Shupp.
\newblock{\em Combinatorial group theory.}
\newblock{Springer-Verlag, Berlin, 1977. 
Ergebnisse der Mathematik und ihrer Grenzgebiete, Band 89.}

\bibitem[Mar95]{Mar:tesi}
R.~Martin.
\newblock {\em Non-uniquely ergodic foliations of thin type, measured currents
  and automorphisms of free groups}.
\newblock PhD thesis, University of California, Los Angeles, 1995.

\bibitem[Whi36]{whi36}
J.~H.~C. Whitehead.
\newblock On equivalent sets of elements in a free group.
\newblock {\em Ann. of Math. (2)}, 37(4):782--800, 1936.

\end{thebibliography}
\end{document}